\documentclass[12pt]{amsart}

\usepackage{cite}

\usepackage{kotex}
\usepackage{lineno,hyperref}
\modulolinenumbers[5]
\usepackage{epstopdf}
\usepackage{pifont}
\usepackage{latexsym}
\usepackage{graphics}
\usepackage{graphicx}
\usepackage{float}
\usepackage[dvips]{xy}
\xyoption{all}
\usepackage{ifpdf}
\usepackage{multirow}
\usepackage{amssymb, amsthm, amsmath}
\usepackage{hhline}
\usepackage{centernot}
\usepackage{verbatim, comment, color}
\usepackage{enumerate}
\usepackage{setspace}
\usepackage{changes}
\usepackage{lipsum}
\usepackage{subcaption}

\newtheorem{theorem}{Theorem} 
\newtheorem{proposition}[theorem]{Proposition}
\newtheorem{lemma}[theorem]{Lemma}
\newtheorem{remark}[theorem]{Remark}

\newtheorem{definition}[theorem]{Definition}
\newtheorem{example}[theorem]{Example}

\newcommand{\ba}{\begin{align}}
\newcommand{\ea}{\end{align}}  

\newcommand{\be}{\begin{equation}}
\newcommand{\ee}{\end{equation}}
\newcommand{\bea}{\begin{eqnarray}}
\newcommand{\eea}{\end{eqnarray}}
\newcommand{\barr}{\begin{array}}
\newcommand{\earr}{\end{array}}
\newcommand{\bn}{\begin{enumerate}}
\newcommand{\en}{\end{enumerate}}
\newcommand{\bi}{\begin{itemize}}
\newcommand{\ei}{\end{itemize}}
\newcommand{\bbbm}{\begin{pmatrix}}
\newcommand{\eeem}{\end{pmatrix}}

\newcommand{\bbS}{{\bf S}}

\newcommand{\cP}{{\cal P}}

\newcommand{\cM}{{\cal M}}

\newcommand{\R}{{\mathbf R}}

\newcommand{\ignore}[1]{}{}
\newcommand{\noin}{\noindent}

\newcommand{\nn}{\nonumber}

\newcommand{\p}{{\partial}}

\newcommand{\dom}{\mathop{\rm dom}}
 
 \newcommand{\Int}{\mathop{\rm int}}

\newcommand{\Prob}{{\mathcal P}}
\newcommand{{\QED}}{{\hfill QED} \smallskip}

\newcommand{\Rn}{{\R^n}}
\newcommand{\spt}{\mathop{\rm spt}}

\newcommand\dist{\hbox{\rm dist}}
\renewcommand{\subset}{\subseteq}
\renewcommand{\supset}{\supseteq}
\renewcommand{\phi}{\varphi}
\newcommand{\cal}{\mathcal}

\newcommand{\pr}{\p_{rel}}

\DeclareMathOperator{\conv}{conv}
  \DeclareMathOperator*{\diam}{diam}

 \DeclareMathOperator*{\argmin}{argmin}
  \DeclareMathOperator*{\argmax}{argmax}
\DeclareMathOperator*{\Var}{{Var}}
\DeclareMathOperator*{\bary}{{\bar x}}
\newcommand{\VVar}{{\textstyle \mathop{\rm Var}_V}}
\numberwithin{equation}{section}
\numberwithin{theorem}{section}
\begin{document}
\title
[Geometrical moment bounds] 
{
{Geometrical bounds for the variance and recentered moments$^*$}
}
 %
%

\thanks{\em TL is grateful for the support of ShanghaiTech University, and in addition, 
to the University of Toronto and its Fields Institute for the Mathematical
Sciences, where parts of this work were performed.  RM  acknowledges partial support of his research by
Natural Sciences and Engineering Research Council of Canada Grant 217006-15. 
The authors are grateful to Guido de Philippis, 
Greg Kuperberg, Tomasz Tkocz, and an anonymous seminar participant 
at Seoul National University for stimulating interactions, and to Hyejung Choi for drawing the figures. 
\copyright 2020 by the authors.
}

\date{\today\\   $^*$The present manuscript is partially based on material
which appeared in an early draft of \cite{LimMcCann20p} and which has been excised from
subsequent versions of that preprint.
}

\author{Tongseok Lim and Robert J. McCann}
\address{Tongseok Lim: Institute of Mathematical Sciences \newline ShanghaiTech University, 393 Middle Huaxia Road, Pudong, Shanghai}
\email{tlim@shanghaitech.edu.cn / tlim0213@outlook.com}
\address{Robert J. McCann: Department of Mathematics \newline University of Toronto, Toronto ON Canada}
\email{mccann@math.toronto.edu}


\begin{abstract}
We bound the variance and other moments of a random vector based on 
the range of its realizations,  
thus generalizing inequalities of Popoviciu (1935) and Bhatia and Davis (2000) 
concerning measures on the line to several dimensions. This is done using convex duality
and (infinite-dimensional) linear programming.

The following consequence of our bounds exhibits symmetry breaking, 
provides a new proof of Jung's theorem (1901),
and turns out to have applications to the aggregation dynamics modelling 
attractive-repulsive interactions:  
among probability measures on $\R^n$ whose support has diameter at most 
$\sqrt{2}$,   we show that the variance around the mean is maximized precisely by those measures which 
assign mass $1/(n+1)$ to each vertex of a standard simplex.
  For $1 \le p <\infty$, the 
$p$-th moment --- optimally centered --- is maximized by the same measures
among those satisfying the diameter constraint.
\end{abstract}
\maketitle

\noindent\emph{MSC2010 Classification 62H05, 49N15, 52A40, 60E15, 90C46}

\section{Introduction}

This article concerns the extension of geometrical variance bounds from one to higher dimensions. 
Let $K \subset \Rn$ be a compact set and
$\Prob(K)$ denote the Borel probability measures supported on $K$.
Let
\begin{eqnarray}
\label{mean}
\bar x(\mu) &:=& \int_\Rn x d\mu(x)
\qquad
\\ \mbox{\rm\ and} \quad
\Var(\mu) &:=&  \int_\Rn |x-\bar x(\mu)|^2 d\mu(x)
\label{variance}
\end{eqnarray}
denote the barycenter (or mean) and the variance of $\mu \in \Prob(K)$.
When  $K:=[\underbar k,\bar k] \subset \R$,  
an inequality due to Bhatia and Davis
\cite{BhatiaDavis00} asserts
\begin{align}\label{BD2}
 \Var(\mu) \le (\bar k-\bar x(\mu))(\bar x(\mu)-\underbar k),
\end{align}
with equality if and only if $\spt \mu \subset \{\underbar k,\bar k\}$.
Optimizing over all possible means yields 
\begin{align}\label{Popoviciu}
 \Var(\mu) \le \frac14 (\bar k-\underbar k)^2,
\end{align}
with equality if and only if $\mu = \frac12(\delta_{\underline k} + \delta_{\bar k})$ ---
a result known since Popoviciu's work \cite{Popoviciu35} on polynomial roots, as explained in
\cite{JensenStyan99}.  We propose to explore higher dimensional, i.e. $n>1$,
generalizations of bounds such as \eqref{BD2}--\eqref{Popoviciu} and their cases of equality. 

In higher dimensions, the shape of the set $K\subset \R^n$ plays a non-trivial role
in the formulation of such a bound.  However, it turns out that the variance maximizing
measures must --- in each case --- be supported on the intersection of $K$ with an
enclosing sphere.  This is the content of our first result,  whose statement requires 
taking the convex envelope of the function
\begin{align}\label{d:phi_K}
\phi_K(x) := 
\left\{
\begin{array}{ccc}
- |x|^2 &  \text{\rm if} & x \in K, \\
+\infty & \text{\rm if} &   x \in \R^n \setminus K.
\end{array}
\right.
\end{align}
Convex envelopes are conveniently expressed using the Legendre transform.

Given a Banach space $Z$ and its dual $Z^*$,  recall the Legendre-Fenchel transform of a function $f:Z \longrightarrow \R \cup \{+\infty\}$
is defined on $Z^*$ by 
\begin{equation}\label{Legendre-Fenchel}
f^*(z^*) := \sup_{z \in Z} z^*(z) - f(z). 
\end{equation}
where $z^*(z)$ denotes the duality pairing.
The double Legendre transform $f^{**}$ is well-known to be the largest lower semicontinuous convex function on $Z^{**}$ whose restriction to $Z$
is dominated by $f$.  
%
Letting $\conv(K)$ denote the smallest closed convex set containing $K$ and $\Int(K)$ the interior of $K$, our 
multidimensional analogs of the Bhatia, Davis \cite{BhatiaDavis00}
and Popoviciu \cite{Popoviciu35} inequalities \eqref{BD2}--\eqref{Popoviciu} are the following:

\begin{theorem}[Enclosing spheres support variance maximizers%
\label{T:Bhatia-Davis}]\ \\
(a) If the measure $\mu \in \Prob(\R^n)$ has barycenter $\bar x(\mu)$ and vanishes outside the compact set 
$K \subset \R^n$, then 
\be\label{Bhatia-Davis}
\Var(\mu) \le  -|\bar x(\mu)|^2 - \phi_{K}^{**}(\bar x(\mu))
\end{equation}
where $\phi_K^{**}$ is defined as in \eqref{d:phi_K}--\eqref{Legendre-Fenchel}.  {If $\bar x(\mu) \in \Int(\conv(K))$,  
then equality holds in \eqref{Bhatia-Davis} if and only if $\mu$ vanishes outside the boundary of 
some closed ball $B$ containing $K$,
i.e. if and only if $\mu[K \cap \p B]=1$.}



(b) Among measures with all barycenters, $\mu$ maximizes variance over $\cP(K) $ if and only if $\mu[K \cap \p B] =1$ where $B$ is the smallest {closed} ball containing $K$, and $\bar x(\mu)$ is the center of $B$. 
Moreover, in this case $\Var(\mu) = R^2$ where $R$ is the radius of $B$.
\end{theorem}

\begin{figure}[b]
\includegraphics[width=0.5\linewidth]{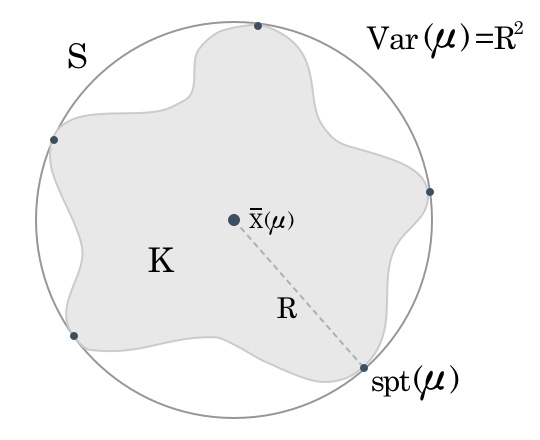}
\caption{Variance maximizer $\mu$ over $\cP(K)$ in Theorem \ref{T:Bhatia-Davis}.}
\end{figure}

{We note that the function $- \phi^{**}_K$ is the concave envelope of $-\phi_K$.} Some refinements and examples include:

\begin{remark}[Specialization to one-dimension]
In the classical context $n=1$ and $K = [\underbar k, \bar k]$,  we recover \eqref{BD2} from \eqref{Bhatia-Davis} by noting
$t\in [0,1]$ and  $x =(1-t)\underbar k +t \bar k$
imply
$$
- \phi_K^{**}((1-t) \underbar k + t \bar k) = (1-t)\underbar k^2 + t \bar k^2
$$
and hence 
$$
- \phi_K^{**}(x) = (\bar k + \underbar k)x - \bar k \underbar k.
$$
\end{remark}

 For our chacterization of equality in \eqref{Bhatia-Davis},
the assumption $\bar x(\mu) \in \Int(\conv(K))$ is in general necessary and cannot be omitted, as the following example indicates.

\begin{example}[Stadium]\label{E:stadium} 
 Taking $K \subset \R^2$ to be the convex hull of two (say unit) balls in $\R^2$ 
and constraining the barycenter $\bar x$ to be (say) the midpoint of one of the flat sides
of $K$ shows the conclusion of Theorem~\ref{T:Bhatia-Davis}(a) need not remain true for all
$\bar x$ in the boundary of $\conv(K)$; the putative enclosing sphere 
degenerates to a halfspace $H \supset K$ in this example, with $\mu[K \cap \p H] =1$
being necessary but not sufficient for equality in \eqref{Bhatia-Davis}.
See the next remark concerning lower dimensional spheres, however.
\end{example}

\begin{remark}[Cases of equality for boundary barycenters]\label{R:boundary barycenters}
 Let $L:= conv(K)$ denote the convex hull of $K \subset \R^n$, i.e.~the smallest closed convex set containing $K$.
 Theorem 2.1.2 of Schneider's book \cite{Schneider93} asserts that each point $x \in L$ belongs to the relative interior of a uniquely determined face $F_x$ of $L$, 
 where a face $F \subset L$ refers to a set containing the endpoints of every segment in $L$ whose midpoint lies in $F$. 
 When $\bar x(\mu) \in \p L$ in Theorem~\ref{T:Bhatia-Davis}, let $j$ denote the dimension of $F_{\bar x(\mu)}$.
 When $j>1$, applying the same theorem to $F_{\bar x(\mu)} \subset \R^j$ instead of $K \subset \R^n$ shows equality holds 
 in \eqref{Bhatia-Davis} if and only if $\mu$ is concentrated on a round sphere $\bbS^{j-1} \subset \R^j$ enclosing $F_{\bar x(\mu)}$.
 When $j=0$ then $\bar x(\mu)$ is an extreme point of $K$,  $\mu$ is a Dirac measure, and \eqref{Bhatia-Davis} becomes an equality.
\end{remark}

\begin{example}[Applications to sample geometries]\label{E:sample geometries}
Computing $\phi_K^{**}$ for special choices of $K$, from
Theorem \ref{T:Bhatia-Davis} and Remark \ref{R:boundary barycenters} we deduce: 
\smallskip

\noin(a) (Ball) If $K= \overline{B_R(0)}$ then $\Var(\mu) \le R^2 - |\bar x(\mu)|^2$,  and equality holds if and only if $\mu$ is supported on $\p K$.
\smallskip

\noin(b) (Ellipse) If $K = \{ (x_1,x_2) \in \R^2 \mid (\frac {x_1} a)^2 + (\frac {x_2} b)^2 \le 1 \}$ with $a >b>0$ 
and   
$\Var(\mu) =  -|\bar x(\mu)|^2 - \phi_{K}^{**}(\bar x(\mu))$ then $\spt \mu$ consists of at most two points.
\smallskip

\noin(c) (Rectangular parallelopiped) If $K = \prod_{i=1}^n [-a_i,a_i]$ is non-empty,  then $\Var(\mu) \le - |\bar x(\mu)|^2 + \sum_{i=1}^n a_i^2$,
and equality holds if and only if $\mu$ is concentrated on the vertices of $K$. 
\smallskip

\noin(d) (Diamond) If $a_1>a_2>0$ and $K = \{ (x_1,x_2) \in \R^2 \mid |\frac{x_1}{a_1}| + |\frac{x_2}{a_2}| \le 1\}$,  then 
$\Var(\mu) \le a_1^2 -\frac{a_1^2 - a_2^2}{a_2} |\bar x_2(\mu)| - |\bar x(\mu)|^2$ and equality holds if and only if 
$\mu$ concentrates at the two vertices of $K$ farthest from the origin, plus at most one of its other two vertices.
\end{example}


Theorem \ref{T:Bhatia-Davis}(b) 
also has analogs for other, possibly anisotropic, measures of the extent to which the mass of $\mu$ 
is concentrated or dispersed. To {illustrate, we give the following definition,  which can be 
contrasted with other generalizations of the variance from the literature, such as those of \cite{PronzatoWynnZhigljavksy17} and its references.}


Let $V:\R^n \to [0,\infty)$ be convex. 
Define
\begin{equation}\label{VVar}
\VVar(\mu) := \inf_{z \in \R^n} \int_\Rn V(x-z) d\mu(x).
\end{equation}
We say that 
$V:\R^n \to [0,\infty)$ is {\em coercive} if its sublevel
sets $V^{-1}([0,\lambda])$ for each $\lambda \ge 0$ are compact.

\begin{remark}[Generalized variances and centered $p$-th moments]
If $V$ is coercive the infimum \eqref{VVar} is attained.
If $V$ is also strictly convex 
and $\VVar(\mu) < \infty$,  then the point $\bar x_V(\mu)$ attaining it is unique,
by displacement convexity \cite{McCann97}.  We can think of $\Var_V(\mu)$ and $\bar x_V(\mu)$  as generalizations of the variance and mean,  which reduce to the classical variance and mean in case 
$V(x)=|x|^2$. When $V(x) = |x|^p$ they reduce to $p$-th moments, but centered on $\bar x_V(\mu)$ 
rather than the classical mean.
\end{remark}

We then generalize Theorem \ref{T:Bhatia-Davis}(b) as follows:

\begin{theorem}[Maximizing generalized variances]
\label{T:generalized variance}
Let $K \subset \R^n$ be compact and $V:\R^n \to [0,\infty)$ be convex and coercive.
Let $\lambda \ge 0$ be the smallest value for which there exists $z \in \Rn$ with
$K+z \subset V^{-1}([0,\lambda])$. 
Then
\begin{equation}\label{max level}
\lambda = \sup_{\mu \in \Prob(K)} \VVar(\mu).
\end{equation} 
Moreover $\mu \in \Prob(K)$ attains this supremum if and only if there exists
\begin{equation}\label{generalized center}
z_* \in \argmin_{z \in \Rn} \int V(x-z) d\mu(x)
\end{equation}
such that $\spt \mu \subset V^{-1}(\lambda)-z_*$.
\end{theorem}

Taking $V(x) = |x|^2$ so that $\VVar=\Var$,  we recognize $\spt \mu \subset V^{-1}(\lambda)-z_*$
as the sphericity condition from Theorem \ref{T:Bhatia-Davis} --- and 
\eqref{generalized center} as the barycenter condition from the same theorem.
More generally,  viewing 
 \eqref{max level} as the value to player 1 of a 
 two-player zero-sum game ---  in which the first player chooses $\mu \in \Prob(K)$
 and the second player, knowing $\mu$, chooses $z \in \Rn$.
 We can interpret \eqref{generalized center}
 as player 2's best response to $\mu$, and $\spt \mu \subset V^{-1}(\lambda) -z_*$ as characterizing
 player 1's best response to $z_*$; together they form the conditions for a saddle-point in the payoff 
 function or  equivalently, for a Nash equilibrium; c.f.~\cite{McCannGuillen13}.

\subsection{Regular simplices maximize moments, given diameter} \ \\

For fixed barycenter $\bar x(\mu)$,  the variance \eqref{variance} is a linear function on the convex
set $\Prob(K)$.  It is thus not surprising that our proof of Theorem~\ref{T:Bhatia-Davis} relies on linear programming duality 
(and convex-concave minimax theory in the case of Theorem \ref{T:generalized variance}). 
A more challenging question is to give sharp bounds on the variance and moments of all measures $\mu$ in a {\em non-convex} set, to which many of the standard techniques in the calculus
of variations \cite{Kawohl85} \cite{McCann97} \cite{BorweinZhu13} no longer apply.

The example which motivated our interest in this problem 
concerns the measures satisfying a diameter bound $\diam[\spt \mu] \le 1$.
Here $\spt \mu$ refers to the smallest closed set containing the full mass of $\mu$.
This question arises as an important special case 
in our work on attractive-repulsive interactions, which addresses the patterns formed by a large collection 
of particles or organisms all preferring to be at distance one from each other \cite{LimMcCann20p}.
We resolve this question below,  by showing among measures $\mu \in \Prob(\R^n)$ 
with $\diam[\spt \mu] \le 1$,  the variance and other moments are maximized precisely when the mass of 
$\mu$ is evenly distributed over the $n+1$ vertices of a regular, unit diameter simplex,  i.e. an equilateral triangle if $n=2$ and a regular tetrahedron if $n=3$.
 
 While it may seem surprising to find this solution
breaks rotational symmetry,  such symmetry breakings undoubtedly bear some responsibility for
the zoo of patterns which emerge from the flocking and swarming models discussed in
\cite{AlbiBalagueCarrilloVonBrecht14} 
\cite{BalagueCarrilloLaurentRaoul13} 
\cite{BertozziKolokolnikovSunUminskyVonBrecht15} 
\cite{BurchardChoksiTopaloglu18} 
\cite{CarrilloFigalliPatacchini17} 
\cite{CarrilloHittmeirVolzoneYao16p} 
\cite{CarrilloHuang17} 
\cite{CarrilloChoiHauray14} 
\cite{ChoksiFetecauTopaloglu15} 
\cite{FellnerRaoul10} 
\cite{FrankLieb18}
\cite{ Lopes19} 
and the references there, of which the present problem represents a limiting case \cite{LimMcCann20p}.  
We were also reminded of the role linear programming 
duality plays in confirming the optimality of sphere packings in certain dimensions 
\cite{OdlyzkoSloane79} \cite{CohnElkies03} \cite{Viazovska17}.

\begin{definition}[Simplices]\label{D:simplices} 
(a) A set $K \subset \R^n$ is called a {\em top-dimensional simplex} if $K$ has non-empty interior and is the convex hull of $n+1$ points $\{x_0, x_1,...,x_n\}$ in $\R^n$.
\smallskip

 \noin(b) A set $K \subset \R^n$ is called a {\rm regular $ k$-simplex} if it is the convex hull of $ k+1$ points $\{x_0, x_1,...,x_{k}\}$ in $\R^n$ satisfying $|x_i-x_j|=d$ for some $d>0$ and all $0 \le i < j \le {k}$. 
 The points $\{x_0, x_1,...,x_{k}\}$ are called {\em vertices} of the simplex. 
 \smallskip
 
\noin  (c) In particular, it is called a {\em unit $k$-simplex} if $d=1$.
\end{definition}

{
\begin{remark}[Regular $n$-simplices $K \subset \Rn$ are top-dimensional] 
\label{R:standard simplex}
A regular $n$-simplex with sidelength $d=\sqrt{2}$ is {linearly} isometric to the following {\em standard simplex} in $\R^{n+1}$
\begin{equation}\label{standard simplex}
\Delta^n := \{ a=\{a_1,...,a_{n+1}\} \in [0,1]^{n+1} \mid \sum_{i=1}^{n+1} a_i =1 \}, 
\end{equation}
which can be verified by simple induction on dimension. 
\end{remark}
}


\begin{figure}[h!]
  \centering
  \begin{subfigure}[b]{0.49\linewidth}
    \includegraphics[width=\linewidth]{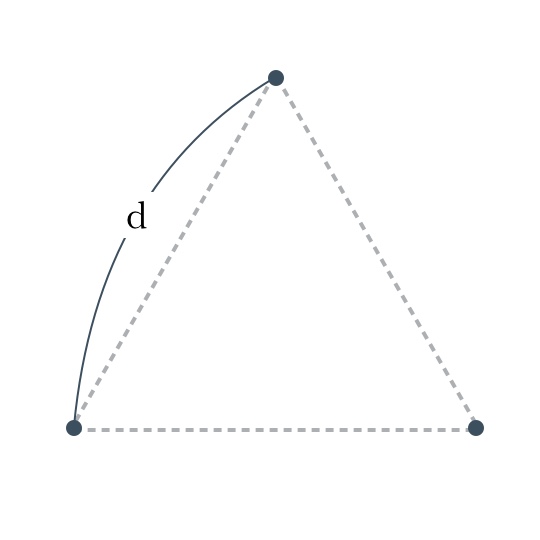}
    \caption{$\spt(\mu)$ in $\R^2$.}
  \end{subfigure}
  \begin{subfigure}[b]{0.49\linewidth}
    \includegraphics[width=\linewidth]{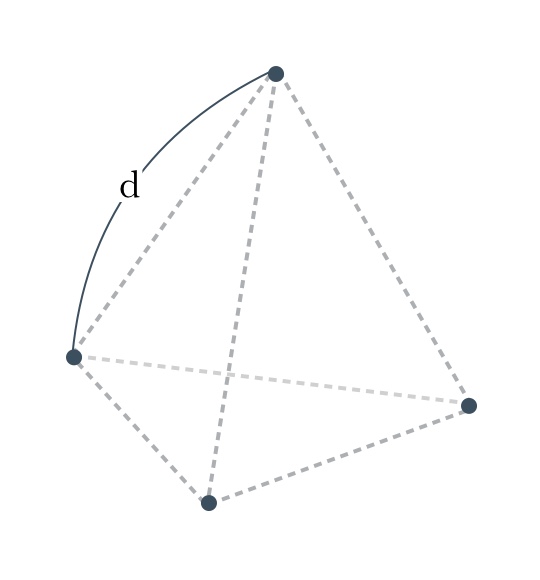}
    \caption{$\spt(\mu)$ in $\R^3$.}
  \end{subfigure}
  \caption{Support of the optimizer $\mu$ in Theorem \ref{T:isodiametric variance bound}.}
  \label{fig:simplex}
\end{figure}

We can now state the following:

\begin{theorem}[Isodiametric variance bounds and cases of equality]
\label{T:isodiametric variance bound}
Let $V(x)=v(|x|)$ with $v:[0,\infty) \longrightarrow [0,\infty)$ convex and increasing.
If the support of a Borel probability measure $\mu$ on $\R^n$ has diameter no greater than $d$, then $\VVar(\mu) \le v (r_n d)$ where $r_n=\sqrt{\frac{n}{2n+2}}$.
Equality holds if and only if 
$\mu$ 
assigns mass $1/(n+1)$ to each vertex of a regular $n$-simplex having diameter $d$. 
\end{theorem}

\begin{example}[Isodiametric bounds on recentered $p$-th moments]
Take $V(x) = |x|^p$ with $p \ge 1$ in Theorem \ref{T:isodiametric variance bound}.
\end{example}

This theorem gives a variational characterization of the unit $n$-simplex.
It can also be viewed as another generalization of Popoviciu's inequality \eqref{Popoviciu}
from $n=1$ to higher dimensions $n>1$. 
 \subsection{Epilog}\ \\
 
After Theorem \ref{T:isodiametric variance bound} 
was announced on the arXiv \cite{LimMcCann20p} (in the special case $V(x)=|x|^2$),
we learned of an isodiametric inequality due to Jung \cite{Jung01}
in which regular simplices also play a crucial role; a modern treatment is given in \cite{DanzerGruenbaumKlee63}:

 \begin{theorem}[Jung] 
\label{T:Jung}
Let $K \subset \R^n$ be compact with $\diam(K)=1$. Then $K$ is contained in a closed ball of radius $r_n=\sqrt{\frac{n}{2n+2}}$. Moreover, unless it lies in some smaller ball, 
$K$ contains the vertices of a unit $n$-simplex. 

\end{theorem}

The constant $r_n$ which appears in these theorems also relates spherical Hausdorff measure
to Hausdorff measure \cite{Federer69}.
Below we shall show how Theorem \ref{T:Jung} follows from our isodiametric variance bound,
thus yielding a new proof of Jung's theorem.  In an appendix to \cite{LimMcCann20p}
we show the converse is also true:  Theorem \ref{T:isodiametric variance bound} 
can be derived from Jung's theorem using elementary geometry. Thus the two theorems are in some sense equivalent.  We are grateful to 
Tomasz Tkocz and an anonymous seminar participant at Seoul National University,
for drawing our attention to Jung's theorem.

\subsection{Plan of the paper:}\ \\

The next section develops the  linear programming and convex duality based proof
of Theorems \ref{T:Bhatia-Davis} and \ref{T:generalized variance}.  
Section \ref{S:variance bound} addresses the non-convex
problem of maximizing moments under a diameter constraint.  It uses induction on dimension to 
prove a geometric lemma which allows us to deduce Theorem \ref{T:isodiametric variance bound},
before closing with 
a new proof of Jung's theorem.

\section{A geometric family of $\infty$-dimensional linear programs}

This section uses linear programming and convex analysis to extend the one-dimensional
inequalities \eqref{BD2}--\eqref{Popoviciu} of Bhatia, Davis and Popoviciu
to higher dimensions, i.e. $n>1$.
Translation invariance allows us to center our measures so that $\bar x(\mu)=0$ without loss of generality.
For each compact $K\subset \R^n$ let
\be\label{zero mean}
\Prob_0(K) := \{ \mu \in \Prob(K) \mid 
\bary(\mu) = 0\}
\end{equation}
denote the set of probability measures on $K$ having vanishing mean.

 Our first goal is to establish the following duality result of Fenchel-Rockafellar type \cite{Rockafellar70}:
 
\begin{proposition}[A strong duality with attainment]\label{P:dual attainment}
If $K \subset \R^n$ is compact then 
\begin{equation}\label{dual attainment}
\sup_{\mu \in \Prob_0(K)} \int_K |x|^2 d\mu(x)
= \inf_{q \in \R^n} \phi_{K}^*(-2q) = - \phi_{K}^{**}(0)
\end{equation}
where  $\phi_{K}^*$ and  $\phi_{K}^{**}$ denote the Legendre transforms \eqref{Legendre-Fenchel} of \eqref{d:phi_K}. 
The supremum is attained if $0 \in L$ and the infimum if $0 \in \Int(L)$, where $L:=\conv(K)$.
A measure $\mu \in \Prob_0(K)$ and point $q \in \R^n$ optimize \eqref{dual attainment}
if and only if $\mu$ vanishes outside  $K \cap \p B_R(q)$ for the smallest sphere $\p B_R(q)$ 
centered at $q$ and enclosing  $K$.
\end{proposition}

Identity \eqref{dual attainment} can be motivated heuristically as follows \cite{McCannGuillen13}.  
Introducing Lagrange multipliers $h$ and $q$ for the mass and barycenter constraints,
 \ba \nn
& \sup_{\mu \in \Prob_0(K)} \int_K |x|^2 d\mu(x) \\
&= \sup_{\mu \in \cM_+(K)}\, \inf_{h \in \R,q \in \R^{n}} h\big(1 - \mu(K) \big) + \int_K \big(|x|^2 +q \cdot x\big)d\mu(x)\nn \\
&\le  \inf_{q \in \R^n, h \in \R}\, \sup_{\mu \in \cM_+(K)} \bigg[ h+  \int_K \big(|x|^2 +q \cdot x -h\big)d\mu(x) \bigg] \nn \\
&= \inf_{q \in \R^n}\, \inf_{h \ge |x|^2+q\cdot x \ \forall x\in K} h \nn \\
&=\inf_{q \in \R^n}\, \sup_{x \in K} \ |x|^2+q\cdot x, \nn
\\ \nn &= -\phi_K^{**}(0)
\end{align}
where $\cM_+(K)$ denotes the set of non-negative Borel measures of finite total mass on $K \subset \R^{n}$.
This inequality can be interpreted as asserting that foreknowledge
of one's opponent's strategy cannot be a disadvantage in a two-player zero-sum game;
it may or may not confer an advantage, depending on the structure of the game.
Statement \eqref{dual attainment} is basically the assertion that the inequality can be replaced with
an equality in our case,
which is a consequence of 
the payoff expression in square brackets having a saddle point  or equivalently, 
of the game having a Nash equilibrium.  
Since the payoff  is bilinear in the variables $\mu$ and $(h,q)$,
this may not be surprising.  Due to lack of compactness however,  a rigorous proof along standard lines requires some machinery. 
%
%
 %
Therefore, recall Theorem~4.4.3 from the book of Borwein and Zhu: 

 \begin{theorem}[Fenchel-Rockafellar duality \cite{BorweinZhu05}] \label{T:BorweinZhu}
Let $A:Z \longrightarrow Y$ be a bounded linear transformation of Banach spaces $Z$ and $Y$,
equipped with functions $f:Z \longrightarrow \R \cup \{+\infty\}$ and $g:Y \longrightarrow \R \cup \{+\infty\}$.
If $g$ is continuous at some point in $A(\dom f)$,  then
$$
 \sup_{y^* \in Y^*} - f^*(A^* y^*) - g^*(-y^*)
= \inf_{z \in Z} f(z) + g(A z),
$$
where $Y^*$ denotes the Banach space dual to $Y$ and $\dom f := f^{-1}(\R)$.
Moreover, the supremum is attained if finite.
\end{theorem}

\noin{\bf Proof of Proposition \ref{P:dual attainment}:} 
Let $Z := \R^{n+1}$ be Euclidean and equip the continuous functions $Y:=C(K)$ on $K$ with the supremum norm,
so that $Z^*=\R^{n+1}$ and $Y^*=\cM(K)$,  the space of signed measures on $K$ normed by total variation.
Take $A(z) = z_0 + \sum_{i=1}^n z_i x_i =: \xi(x) \in Y$ so that $A^*\mu = \int_K (1,x) d\mu(x)$ gives the mass and barycenter of $\mu \in \cM(K)$.
Set $f(z_0,\ldots,z_n) := z_0$ so that
$$f^*(z^*) = 
\left\{ 
\begin{array}{cl} 
0 & {\rm if}\ z^* = (1,0,\ldots,0), \\
\infty & {\rm else}.
\end{array}
\right.
$$
Also set
$$g(\xi) :=
\left\{ 
\begin{array}{cl} 
0 & {\rm if}\ \xi(x) \ge |x|^2 \quad \forall x \in K, \\
\infty & {\rm else},
\end{array}
\right.
$$
so that 
$$g^*(\mu) =
\left\{ 
\begin{array}{cl} 
\int_K |x|^2 d\mu(x) & {\rm if}\ \mu \le 0, \\
\infty & {\rm else}.
\end{array}
\right.
$$
Inserting these choices into Theorem \ref{T:BorweinZhu} 
yields  \eqref{dual attainment},
noting the definitions \eqref{d:phi_K}--\eqref{Legendre-Fenchel} of $\phi_K^*$.
If $0 \in L:=\conv(K)$ then $\Prob_0(K)$ is non-empty and the supremum is bounded above and below (by the infimum and zero) 
hence attained (also by Theorem \ref{T:BorweinZhu}).

 If $0 \in \Int(L)$ then $rB \subset L$ for $r>0$ sufficiently small,
where $B:=B_1(0)$ is the centered unit ball.
Then $\phi_L \le \phi_{rB}$ hence $\phi_L^*(q) \ge \phi^*_{rB}(q) = r|q| + r^2$ grows without bound as $|q| \to \infty$.
Being lower semicontinuous,  $\phi_L^*$ then attains its minimum.  On the other hand,  the concavity of $x \mapsto -|x|^2$
implies $\phi_L^*=\phi_K^*$, as we now argue. Indeed $\phi_L^* \ge \phi_K^*$ follows directly from $K \subset L$ and $\phi_K \ge \phi_L$.
Conversely,  given any affine function $a$ on $\R^n$ dominated by $\phi_K$,  we find $a \le \phi_L$ also,
since $\phi_K=\phi_L$ outside $L\setminus K$, and each $x \in L \setminus K$ can be approximated by convex combinations
$x^{j} = \sum_{i=1}^j t_i^{j} k_i^{j}$ of points $k_i^{j} \in K$ with $t_i^j \ge 0$ and $\sum_{i=1}^j t_i^{j}=1$, so
\begin{align*}
\phi_L(x^{j})
&= -|x^{j}|^2 
\\ & \ge -  \sum_{i=1}^j t_i^{j} |k_i^{j}|^2
\\ & =  \sum_{i=1}^j  t_i^{j} \phi_K(k_i^{j})
\\ & \ge a(x^{j}),
\end{align*}
and the limit $x=\lim_{j\to \infty} x^j$ yields $\phi_L \ge a$ as desired.
Since $\phi_K^{**}$ is the supremum of such affine functions $a$, we conclude $\phi_L^{**} \ge \phi_{K}^{**}$,
which implies $\phi_K^* \ge \phi_L^*$ hence  $\phi_K^* = \phi_L^*$.

To characterize the optimizers,  let $\mu \in \Prob_0(K)$ and $q \in \R^n$.  Then
\begin{align}\nn
\int_K |x|^2 d\mu(x)
&=\int_K (|x-q|^2 - |q|^2) d\mu(x)  
\\ &\le \max_{x \in K} |x-q|^2 - |q|^2
\label{nc}
\\ &= \phi_K^*(-2q)
\nn\\ &=: R^2 - |q|^2
\nn
\end{align}
and equality holds if and only if $\mu$ vanishes outside the set 
$$\argmax_{x \in K} |x-q|^2 = K \cap \p B_R(q);$$
here $R$ is the smallest radius for which $K \subset \overline{B_R(q)}$. 
On the other hand,  $\mu \in \Prob_0(K)$ and $q \in \R^n$ optimize \eqref{dual attainment}
if and only if equality holds in \eqref{nc}, so the proposition is established.
\QED

Expression \eqref{dual attainment} 
is particularly convenient for selecting the translation of $K$ which maximizes 
 the value of the linear program using the following lemma:

\begin{lemma}[Optimal translation of a domain relative to the origin] 
\label{L:optimal translation} 
{For compact $K \subset \R^n$, we have 
$\phi^{**}_{K-w}(x) =  (|x+w|^2 -|x|^2) + \phi^{**}_K(x+w)$.} In particular,
 $\phi^{**}_{K}(0) \le \phi^{**}_{K-w}(0)$ for all {translations} $w \in \R^n$ if and only if $\phi^{**}_K$ 
attains its minimum at the origin.
\end{lemma}

\noin{\bf Proof.}  {The Legendre-Fenchel transform \eqref{Legendre-Fenchel}, applied to $\phi_K$, yields }
\begin{eqnarray*}
 \phi^*_{K-w} (y) 
 &=& |w|^2 -w \cdot y + \phi^*_K(y-2 w) \qquad{\rm and}
 \\
\phi^{**}_{K-w}(x) 
&=&  |w|^2 +2 w\cdot x + \phi_K^{**}(x+w),
\end{eqnarray*}
hence
\begin{align}\label{identity}
\phi^{**}_{K-w}(0) = |w|^2 + \phi_K^{**}(w).
\end{align}
Recall that a convex function $f$ on a Banach space $Z$ attains its minimum at $x$ if and only if $0 \in \p f(x)$, where
\begin{equation}\label{subgradient}
\p f(x) := \{z^* \in Z^* \mid f(z)  \ge f(x) + z^*(z-x) \quad \forall z \in Z\}.
\end{equation}
The formula above shows $f(w) := \phi^{**}_{K-w}(0)$ to be a strictly convex function of $w$ with $\p f(0) = \p \phi_K^{**}(0)$,
so 
{$\phi_{K-w}^{**}(0)$ attains its minimum at $w=0$ if and only if $\phi_{K}^{**}(w)$ does as well. }
 \QED \\

\noin{\bf Proof of Theorem \ref{T:Bhatia-Davis}:} (a)
For a compact set $K\subset \Rn$ with $w \in \R^n$,
Lemma \ref{L:optimal translation} yields
$
\phi^{**}_{K-w}(0) = |w|^2 + \phi_K^{**}(w).
$ 
In \eqref{dual attainment} 
this gives
\begin{align}\label{BD3}
\sup_{\nu \in \cP_0(K-w)} \int |x|^2 d\nu = - \phi^{**}_{K-w}(0) =  -|w|^2 - \phi_{K}^{**}(w).
\end{align}
Letting $\mu$ denote the translation of $\nu$ by $w:=\bar x(\mu)$ yields \eqref{Bhatia-Davis}.
If $w \in \Int(\conv(K))$, Proposition \ref{P:dual attainment} states that $\nu \in \cP_0(K-w)$ attains the supremum if and only
if $\nu$ vanishes outside $K \cap \p B_R(q)$ with $K \subset \overline{B_R(q)}$ for some $q \in \R^n$ and $R>0$.
In other words, if $\bar x(\mu) \in \Int(\conv(K))$, then \eqref{Bhatia-Davis} becomes an equality if and only if
$\mu$ is supported on the boundary of a closed ball containing $K$. 

(b) Can be proved using Lemma \ref{L:optimal translation} as in \cite{LimMcCann20p}. Alternately (b)
 follows from the choice $V(x)=|x|^2$ in Theorem \ref{T:generalized variance}, whose proof appears just below.
\QED


 \noin{\bf Proof of Theorem \ref{T:generalized variance}.} 
Recall from e.g. \cite{McCannGuillen13} that
\begin{align*}
 \sup_{\mu \in \cP(K)} \VVar(\mu) &= \sup_{\mu \in \cP(K)}  \inf_{z \in \R^n} \int V(x-z) d\mu(x) \\
&\le \inf_{z \in \Rn} \sup_{\mu \in \cP(K)}   \int V(x-z) d\mu(x) \\
&= \min_{z \in \Rn} \max_{x \in K} V(x-z) 
\\ &=\lambda.
\end{align*}
Combining compactness of $K$ with coercivity and continuity of the convex function $V$ allows us to replace 
$\R^n$ with a sufficiently large closed ball $\overline{B_R(0)}$ without affecting the values of either infimum;
the infima are therefore attained, and the
inequality above becomes an equality according to convex-concave minimax theory, 
e.g.~\cite[Theorem 45.8]{Strasser85}. 

From the definition of $\lambda$,  there exists $z_*$ such that $K-z_* \subset V^{-1}([0,\lambda])$.
Thus $\mu \in \cP(K)$ satisfies
$$
\inf_{z \in \Rn} \int_K V(x-z) d\mu(z) \le \int_K V(x-z_*) d\mu(x) \le \lambda
$$
with the first inequality being saturated if and only if \eqref{generalized center} holds,
and the second inequality being saturated if and only if $V(x-z_*) = \lambda$ on $\spt \mu$.
In light of \eqref{max level},
these two conditions are necessary and sufficient to ensure that $\mu$ is a maximizer.
\QED

 \section{Isodiametric variance and $p$-th moment bounds}
 
{
\label{S:variance bound}

This section establishes our isodiametric variance bound and cases of equality: Theorem \ref{T:isodiametric variance bound}.
Let us briefly outline the strategy of our proof.   Fix $V(x)=v(|x|)$ radially symmetric, convex and increasing.
For each compact set $K\subset \R^n$ of unit diameter,  
Theorem \ref{T:generalized variance} asserts
 (i) that the maximizer of $ \VVar(\mu)$
on $\Prob(K)$ vanishes outside the smallest sphere enclosing $K$ and 
(ii) the center of this sphere attains the infimum \eqref{VVar} defining $\Var_V(\mu)$.
We may, without loss of generality assume that $K$ has been translated
so that this sphere is centered on the origin. 
We shall now show the radius of this sphere cannot exceed the radius $r_n := \sqrt{\frac n{2n+2}}$ 
of the unit $n$-simplex. To do so we use 
an induction on dimension, 
which is based on the idea that if the centered sphere is too large,  no measure whose support has unit diameter can have its center of mass at the origin.  More precisely, we show the following elementary yet crucial geometric proposition which 
characterizes the unit simplex.

 
 \begin{proposition}[Tension between diameter and center-of-mass constraints] 
 \label{P:tension}
(a) If $K \subset \p B_r(0)$ 
is a subset of the radius $r >r_n := \sqrt{\frac{n}{2n+2}}$
centered sphere in $\Rn$ 
and ${\rm diam}(K) \le 1$, then $0 \notin \conv(K)$.   
\vspace{1mm}
 
\noin (b) If $K$ is a subset of the centered sphere  in $\R^n$ of radius $r_n$, ${\rm diam}(K) \le 1$ and $0 \in \conv(K)$, then $K$ is the set of vertices of a unit $n$-simplex.
\end{proposition}

\noin{\bf Proof of Proposition \ref{P:tension}.} {\bf (a)} The proposition is trivial to verify when $n=1$. To derive a contradiction, suppose the proposition holds in $\R^{n-1}$ 
but fails in $\R^n$. 
Then there exists a centered sphere $S$ of radius $r$ with $r > r_n$, and $K \subset S$ with ${\rm diam} (K) \le 1$ and $0 \in \conv(K)$. We can find $n+1$ points in $K$, say $X:= \{x_0,x_1,...,x_n\} \subset K$, such that $0 \in \conv(X)$.  If the origin lies on the boundary of $\conv(X)$, then after 
intersecting the problem with a hyperplane supporting $\conv(X)$ at $0$, 
the inductive hypothesis yields the desired contradiction using $r_{n-1}<r_n$. We may therefore assume
$0 \in  {\rm int}\conv(X)$, so that $\conv(X)$ is a top-dimensional simplex in $\R^n$. 

Without loss of generality, let $x_0=r\hat e_1=(r,0,...,0)$.  Define
\[  U:=\{x \in S \ | \ |x - x_0| \le 1\}.
\]
Then $\pr U := \{x \in S \ | \ |x - x_0| = 1\}$ is a $(n-2)$-dimensional sphere of radius $r'$ and center {$a=a_1\hat e_1$ for some $r' >0$ and  $a_1 \in \R$.} 
Since $0 \in  {\rm int}\conv(X)$ implies $0 \in  {\rm int}\conv(U)$, we see that $a_1 <0$.
And $r > r_n$ implies $r' > r_{n-1}$, as $r' = r_{n-1}$ precisely when $r = r_n$. 
Now consider the unique hyperplane $H$ which contains the $(n-1)$-{simplex} with vertices $X'=\{x_1,...,x_n\} \subset X$. 
Let $L$ be the one-dimensional subspace spanned by $\hat e_1$. Then $H \cap L \neq \emptyset$ since $0 \in  {\rm int}\conv(X)$. 
Let $b=b_1\hat e_1:=H \cap L$. Then $a \le b_1$ since $X' \subset U$, and $b_1 < 0$ since $0 \in  {\rm int}\conv(X)$. 
Now define the disk $D:=\conv(H \cap S)$ whose (relative) boundary is 
the $(n-2)$-dimensional sphere $\pr D := H \cap S$. Note that $b \in D$ and $X' \subset \pr D$.  Define 
\[ d := \dist (b, \pr D).
\]
Notice that the facts $a_1 \le b_1 < 0$ and $\pr D \subset U$ imply $d \ge r'$, hence $d > r_{n-1}$ {(see Figure \ref{fig:circle})}. 



The desired contradiction (and proposition) will follow if we show that $b \notin \conv (X')$, as this will  imply $0 \notin \conv (X)$. 
 To achieve this,  suppose on the contrary $b \in \conv (X')$.  
Let $D'$ be the $(n-1)$-dimensional closed ball in $H$ of center $b$ and radius $d$, and let $\pr D'$ be its boundary sphere.
Note that $b \in \conv (X') \cap D'$. Since none of the extreme points of {$\conv(X')$} lie in the interior of $D'$,
it follows the extreme points of $\conv(X') \cap D'$ all lie on the boundary sphere $\pr D'$.  Setting {$K' := \conv(X') \cap \pr D'$,}
the Krein-Milman theorem implies $b \in \conv(K')$.  But this contradicts the inductive hypothesis,  which asserts that the center $b$
of a sphere $S':= \pr D'$ of radius $d > r_{n-1}$ cannot lie in the convex hull of any subset $K' \subset S'$ whose diameter is bounded by one. 
\begin{figure}[h!]
  \centering
  \begin{subfigure}[b]{0.49\linewidth}
    \includegraphics[width=\linewidth]{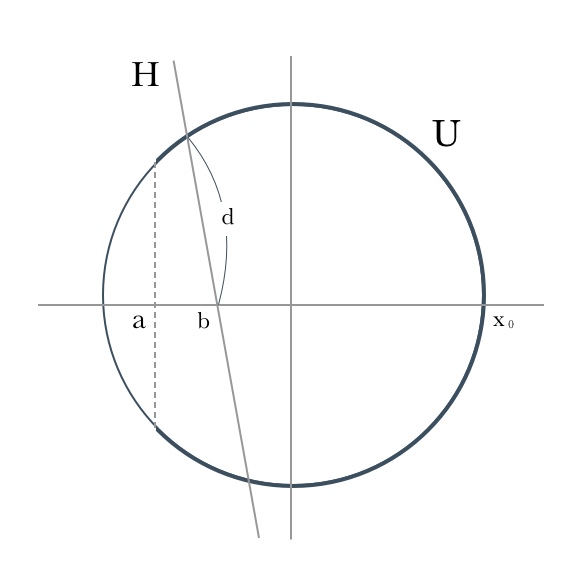}
    \caption{$a, b \in \R^n$ and $d >0$.}
  \end{subfigure} 
  \begin{subfigure}[b]{0.49\linewidth}
    \includegraphics[width=\linewidth]{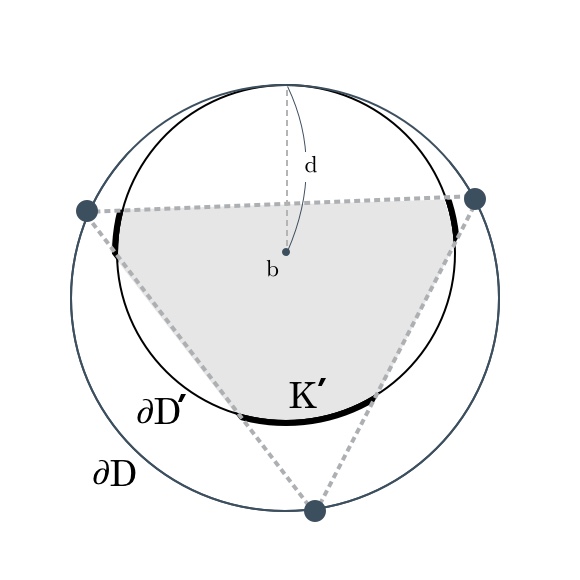}
    \caption{$b \in \conv(K')$.}
  \end{subfigure}    \caption{$b \in \conv(K')$ yields a contradiction.}
  \label{fig:circle}
\end{figure}

\noin {\bf (b)} We proceed as in part (a). Suppose the proposition holds in $\R^{n-1}$. Let $S$ be the centered sphere of radius $r_n$ in $\R^n$, and let $K \subset S$ be such that ${\rm diam} (K) \le 1$ and $0 \in \conv(K)$.  As before we can find a subset $X$ of $K$, the vertices of a $n$-simplex  with $0 \in  \conv(X)$, and in fact $0 \in  {\rm int}\conv(X)$ {by part (a)}. 
Note that the sphere $\pr U$ now has radius $r_{n-1}$. Again consider the hyperplane $H$ spanned by $X'$, and observe that $b=b_1 \hat e_1 \in \conv (X')$ since $0 \in \conv (X)$. Now if $a_1<b_1$, then as before we have $d > r_{n-1}$. This yields a contradiction by part (a) and the last part of its proof. We conclude that $a_1=b_1$, and this implies that $H$ is the hyperplane containing $b$ and having $x_0=r_n \hat e_1$ as its normal. Then {$X' \subset H \cap S = \pr U$,} and  the induction hypothesis implies that $X'$ must form vertices of a  unit $(n-1)$-simplex. Hence $X$ forms vertices of a unit $n$-simplex, inscribed in the sphere $S=\p B_{r_n}(0)$.

It remains to show that $K=X$.  
Since $\conv(X)$ is an intersection of $n+1$ closed halfspaces and $X=\conv(X) \cap S$,  any point $x' \in K \setminus X$ lies outside 
at least one of these halfspaces.  Without loss of generality, we may suppose it lies in the halfspace $H_a := \{ x \in \R^n \mid x \cdot \hat e_1 < a_1\}$. 
But this means $x' \in S \setminus U$, yielding $|x'-x_0| >1$, which contradicts the assumption $\diam(K) \le 1$. 
\QED


We are now in a position to prove Theorem \ref{T:isodiametric variance bound} by characterizing variance maximizing measures under a diameter constraint.\\

\noin{\bf Proof of Theorem \ref{T:isodiametric variance bound}}: 
Set $V(x)=v(|x|) \ge 0$ with $v$ convex and increasing, and
fix a compact set $K \subset \R^n$ with diameter no {greater} than $1$,
and let $\mu \in \Prob(K)$ be the probability measure on $K$ 
which maximizes $\Var_V.$ 
Such a measure exists,
by the weak-$*$ compactness of $\Prob(K)$ in the Banach space $\cM(K)$ dual to $(C(K),\| \cdot \|_\infty)$
(or by Proposition \ref{P:dual attainment} in case $V(x) = |x|^2$).
We may assume $K$ 
has been translated so that the origin $z_* =0$ satisfies \eqref{generalized center}.
In this case we claim $0 \in \conv(\spt \mu)$.  If not, letting $0 \ne z$ be the point of $\conv(\spt \mu)$ closest
to the origin,  say $z= (r,0,\ldots,0)$,
we find each point $x \in \conv(\spt \mu)$ lies in the halfspace to the right of $z$, hence is strictly closer to $z$ than to $0$,
contradicting \eqref{generalized center}.
Theorem~\ref{T:generalized variance} asserts $\mu$ vanishes outside the smallest sphere
$\overline{B_R(0)}$ enclosing $K$,  so that $\VVar(\mu) = v(R)$.
On the other hand, $\spt \mu \subset \p B_R(0)$  has diameter at most one and contains 
the origin in its convex hull.
Proposition~\ref{P:tension} therefore asserts that $R \le r_n$ 
and that when equality holds {$\spt \mu$} coincides
with the vertices of a unit $n$-simplex.   Note that the uniform measure $\hat \mu$ on the vertices of this simplex has center of mass
at the origin and $\VVar(\hat \mu) = v(r_n)$.  Remark \ref{R:equidistribution} below shows
no other measure on the vertices of the simplex has center of mass at the origin.
 If $R<r_n$ we conclude $\VVar(\mu)<\VVar(\hat \mu)$,
while if $R=r_n$ we conclude $\mu=\hat \mu$.  Thus for the given diameter $d=1$ of support, we have identified the maximum of $ \VVar(\cdot)$ 
and the measures which attain it uniquely {(up to translations and rotations)}.
\QED

\begin{remark}[Equidistribution over the simplex vertices]\label{R:equidistribution}
Since the vertices of the standard simplex 
\eqref{standard simplex} form a basis for $\R^{n+1}$,
each point inside the simplex can be uniquely expressed as a convex combination of its vertices.
Thus {among measures on the vertices of 
the simplex, only the uniform measure has its barycenter at the 
point} $\frac1{n+1}(1,\ldots,1)$. 
\end{remark}
 
 \subsection{A new proof of Jung's theorem}\ \\
 
Let conclude by showing how Jung's theorem \cite{Jung01} follows from the results just derived:
 \\
 
\noin{\bf Proof of Theorem \ref{T:Jung} using Theorems \ref{T:isodiametric variance bound} and  \ref{T:Bhatia-Davis}(b).} Let $K \subset \R^n$ 
be compact with $\diam(K) \le 1$. Theorem \ref{T:isodiametric variance bound} asserts that any $\mu \in \cP(K)$ satisfies $\Var(\mu) \le r_n^2$, and Theorem \ref{T:Bhatia-Davis}(b) then implies that $K$ can be contained in a closed ball of radius at most $ r_n$. Now suppose $K$ does not lie in a ball with radius strictly smaller than $r_n$. Then Theorem~\ref{T:Bhatia-Davis}(b) provides $\mu \in \cP(K)$ with $\Var(\mu) = r_n^2$ and Theorem \ref{T:isodiametric variance bound} then implies that $\spt(\mu)$ contains the vertices of a unit $n$-simplex. \QED

Conversely,  an appendix to our companion work \cite{LimMcCann20p}
{shows how Jung's theorem can be used to prove 
Theorem \ref{T:isodiametric variance bound}  --- at least for $V(x)=|x|^2$, but
the proof there adapts easily to other radially symmetric, convex increasing choices of  
$V(x)=v(|x|)$.}
 



\end{document}